\documentclass[12pt]{article}
\usepackage{amssymb,amsmath,amscd}

\newtheorem{proposition}{Proposition}
\newtheorem{lemma}[proposition]{Lemma}

\def\A{{\cal A}}
\def\Co{\mathbb{C}} 
\def\fidi{\hskip5pt \vrule height4pt width4pt depth0pt \par}
\def\id{{\rm id}}
\def\Na{\mathbb{N}}  
\def\Z{\mathbb{Z}} 
\def\R{{\cal R}}
\def\Re{\mathbb{R}}  
\def\sphere{\mathbb{S}}
\def\sph4{\Sigma_q^4}
\def\ot{\otimes}
\def\lra{\longrightarrow}
\def\dr{\mbox{$\Delta_{R}$}}
\def\o{\sp{[1]}}
\def\t{\sp{[2]}}

\newcommand{\nc}[2]{\newcommand{#1}{#2}}
\nc{\inc}{\mbox{$\,\subseteq\;$}}

\begin{document}

\title{Bijectivity of the canonical map for the noncommutative instanton bundle}

\author{F.Bonechi${}^{1,2}$, N.Ciccoli${}^3$, L.D\c abrowski${}^{4}$,
M.Tarlini${}^{1,2}$}

\date{July 14, 2003}

\maketitle

\centerline{{\small  ${ }^1$ INFN Sezione di Firenze}}
\centerline{{\small ${ }^2$ Dipartimento di Fisica, Universit\`a
di Firenze, Italy. }}
\centerline{{\small ${ }^3$ Dipartimento di
Matematica, Universit\`a di Perugia, Italy. }}
\centerline{{\small ${ }^4$ Scuola Internazionale Superiore di Studi Avanzati,
Trieste, Italy. }}
\centerline{{\small
e-mail: bonechi@fi.infn.it, ciccoli@dipmat.unipg.it,
dabrow@sissa.it, tarlini@fi.infn.it}}

\vskip2cm
\begin{abstract}
\noindent
{It is shown that the quantum instanton bundle introduced
in \cite{BCT1} has a bijective canonical map and is, therefore, a coalgebra
Galois extension.}
\end{abstract}

\medskip
\vskip1cm
\noindent
{\bf Math.Subj.Classification}: 16W30, 19D55, 58B32, 58B34, 81R50


\thispagestyle{empty}

\bigskip
\bigskip

\newpage
\section{Introduction}

Since the beginning of the theory of quantum groups lot of efforts
have been devoted to develop a full quantum analogue of the notion of
principal bundle. For that purpose the key feature turns out to be
a characterization of freeness and transitivity on the fiber of the
group action through the bijectivity of the so--called canonical map.
In \cite{Sch} the relevant role of this map in the context of Hopf
algebras was stressed and a theory of algebraic quantum principal bundles
with a Hopf algebra as structure group was built up. A notion related
to such a Hopf-Galois property but on the level of differential calculi
and principal connections was partially used in \cite{BM0,D1}. Later
on several examples of quantum principal bundles have been shown to
be indeed Hopf--Galois extensions and thus to fully deserve their
name, see e.g. \cite{D} for a list.

However many interesting examples of deformations of principal bundles
coming from quantum groups could not be expressed in that framework.
For instance a quantum homogeneous space is a Hopf--Galois extension
only if it is a quotient by a quantum subgroup, which is a very rare case.
It is therefore necessary to consider a more general object in place of
the structure group. In \cite{Sch} bundles with a \emph{coalgebra structure
group} obtained as quotient of Hopf algebras by a coideal and right ideal
were first considered. Later on it was shown in \cite{Brz} that all Podle\'s
spheres can be obtained as quotients by a \emph{coalgebra subgroup}
of $SU_q(2)$.
Such results were further developed in \cite{Brz2,BH1,BM} to the idea of
coalgebra--Galois extensions and, more recently, in \cite{BH2} into a
possible definition of coalgebra principal bundles.
It has to be stressed that the bijectivity of the canonical map
retains a crucial role in the theory.

Approximately at the same time a \emph{semiclassical} ({\it i.e.}
Poisson) interpretation for this approach was given in \cite{Ci}
and \cite{BCGST}, in terms of coisotropic subgroups of a
Poisson--Lie group. The semiclassical point of view turned out to
be very useful as it allowed in \cite{BCT1} to describe a quantum group
version of the $SU(2)$--principal bundle $\sphere^7\to\sphere^4$.
In \cite{BCT2} a description of this bundle in terms of quantized
enveloping algebras is given.

It is then quite natural to ask whether this construction can be
considered an algebraic quantum principal bundle, the key point being
a verification of the bijectivity of the canonical map, as explained.
In this paper we will prove that property.
Let us remark that the quantum instanton bundle is far more difficult to
deal with than any other known example arising from quantum groups,
for at least three different reasons: the base space is not a quantum
homogeneous space (just a double coset of $U_q(4)$), the structure
group is only a coalgebra and, lastly, such coalgebra corresponds
to a non abelian group.

Once the Galois property is proven the next natural step to
understand the geometry of this quantum principal bundle is
to analyze the quantum vector bundles associated to
corepresentations of the structure coalgebra. The bijectivity of the
canonical map proven in this paper together with the results contained
in \cite{BH2} and in the work in preparation \cite{schsch} allow to
conclude that they are finitely generated and projective modules. In the
last section we will explain this point.

Our results and the explicit expression of the $K$--homology
generators given in \cite{BCT1} are all the needed ingredients to apply the
character formula in \cite{BH2} to compute the charges of these bundles,
in analogy to the computations carried out for line bundles on Podle\'s
sphere, both in the Hopf--Galois \cite{Haj} and in the coalgebra--Galois
case \cite{HMS}.


\section{The quantum four-sphere $\Sigma^4_q$}
In this section we recall the basic facts necessary to obtain the
four-sphere $\Sigma^4_q$ (\cite{BCT1}). The algebra of polynomial
functions $\A(U_q(4))$ is generated by $\{t_{ij}\}_{ij=1}^4$,
$D_q^{-1}$ and the following relations:
\begin{equation}
\label{cr}
\begin{array}{rcll}
t_{ik} t_{jk} = q \ t_{jk} t_{ik}\;,
& {} & t_{ki} t_{kj} = q \ t_{kj} t_{ki}\;,  &i<j \cr
t_{i\ell} t_{jk} &=&  t_{jk} t_{i\ell}\;, &i<j,\,k<\ell\cr
t_{ik} t_{j\ell} - t_{j\ell} t_{ik}
& = & (q-q^{-1}) t_{jk} t_{i\ell}\;, &i<j,\, k<\ell\cr
D_qD_q^{-1} &=& D_q^{-1} D_q = 1\;,&{}
\end{array}
\end{equation}
where $D_q =\sum_{\sigma\in
P_4}(-q)^{\ell(\sigma)}t_{\sigma(1)1}\ldots t_{\sigma(4)4}$ is the
quantum determinant and $P_4$ is the group of 4--permu\-tations.
It is easy to see that $D_q$ is central. The coproduct is
\begin{equation}
\label{has}
\Delta(t_{ij})=\sum_k t_{ik}\otimes t_{kj}\;,\ \
\Delta(D_q)=D_q\otimes D_q\; .
\end{equation}
For the usual definition of the antipode $S$ see \cite{Schmdg};
the compact real structure forces to choose $q\in\Re$ and it reads
\begin{equation} \label{rs}
t^*_{ij} = S(t_{ji}), ~~~~~D_q^* = D_q^{-1}\;.
\end{equation}

In the following we will denote $\kappa=*\circ S$.

The algebra of polynomial functions $\A(\sphere^7_q)$ on the quantum
seven-sphere (see \cite{VS}) is generated as a $*$-algebra by
$\{z_i = t_{4i}\}_{i=1}^4$, which verify the following relations
\begin{equation}
\label{s7}
\begin{array}{c}
z_i z_j\ =\ q z_j z_i  ~~(i<j) ~,~~~~~~~~~ z_j^*z_i = q z_i z_j^*
~~(i\neq j)~, \cr z^*_k z_k = z_k z^*_k + (1-q^2) \sum_{j<k} z_j
z_j^* ~,~~~~~~~~~\sum_{k=1}^4 z_k z^*_k = 1\;.
\end{array}
\end{equation}
The $\A(U_q(4))$-coaction on $\sphere^7_q$ reads
\begin{equation}
\label{coU}
\Delta(z_i)=\sum_j z_j\otimes t_{ji}~, ~~\Delta(z_i^*)=\sum_j z_j^*\otimes t_{ji}^*\ .
\end{equation}

Define $\R=R\ \A(U_q(4))$, where
\begin{eqnarray*}\label{ideal1}
R&=&{\rm Span}\{t_{13},\,t_{31},\,t_{14},\,t_{41},\,t_{24},\,t_{42},
\,t_{23},\,t_{32},\,t_{11}-t_{44},\,t_{12}+t_{43},\cr
&&\qquad\ \ t_{21}+t_{34},\,t_{22}-t_{33},\,t_{11}t_{22}-q\ t_{12}t_{21}-1\}\ .
\end{eqnarray*}
It is easy to verify that $\R$ is a $\kappa$--invariant, right ideal,
two sided coideal. In the sequel we shall denote $C := \A(U_q(4))/\R$ the
quotient and $r:\A(U_q(4))\to C$ the canonical projection.

By construction $C$ is a coalgebra, a right $\A(U_q(4))$--module
and it inherits an involutive, antilinear map $\kappa_C$.
In \cite{BCT1} it has been shown that $C$ is isomorphic to$\A(SU_q(2))$
as coalgebras, and that the isomorphism intertwines $\kappa_C$ with $*\circ S$
on $\A(SU_q(2))$. Using this isomorphism we could transfer to $C$ the algebra
structure of $\A(SU_q(2))$ (and thus of a Hopf-algebra) but the projection map
$r:\A(U_q(4))\to C$ would not be a homomorphism,
{\it e.g.} $r(t_{11} t_{43})\not= r(t_{11}) r(t_{43})$.
In the rest of the paper we actually assume this point of view,
and identify $C$ with $\A(SU_q(2))$. Although $r$ doesn't respect the algebra
structure $C$ has a right $\A(U_q(4))$--module structure, that will be important
for us, defined by
\begin{equation}
\label{Ums} r(t)\, \cdot t' := r(tt')\ .
\end{equation}
We introduce in $C$ the usual $\Co$-linear basis
$\langle r^{k,m,n}\ |\ k\in\Z;\
m,n\in\Na\rangle$, where
\begin{equation}
\label{su2basis}
r^{k,m,n} = \left\{
\begin{array}{ll}
r(t_{11}^k t_{12}^m t_{21}^{n})&\quad k\geq0\;, \\
r(t_{12}^m t_{21}^n t_{22}^{-k})&\quad k <0\;.
\end{array}
\right.
\end{equation}

\medskip
In the following Lemma we collect some useful relations concerning
the $\A(U_q(4))$--module structure of $\A(SU_q(2))$.

\begin{lemma}
The following relations are valid for $k\in\Z$ and
$m,n\in\Na$:
\begin{eqnarray}
\label{modulo_relazioni} r^{k,m,n}\cdot\left(\begin{array}{cc}
                                t_{11}^*&t_{12}^*\cr
                                t_{21}^*&t_{22}^*
                                \end{array}\right)&=&
r^{k,m,n}\cdot\left(\begin{array}{cc}
                                t_{22}&-qt_{21}\cr
                                -q^{-1}t_{12}&t_{11}
                                \end{array}\right)          \cr
r^{k,m,n}\cdot\left(\begin{array}{cc}
                                t_{33}^*&t_{34}^*\cr
                                t_{43}^*&t_{44}^*
                                \end{array}\right)&=&
r^{k,m,n}\cdot\left(\begin{array}{cc}
                                q^{m+n}t_{11}&q^{1-k}t_{12}\cr
                                q^{-(1+k)}t_{21}&q^{-(m+n)}t_{22}
                                \end{array}\right)          \cr
\nonumber
& &-\theta(k)q(1-q^{-2k})r^{(k-1),m,n}\cdot\left(\begin{array}{cc}
                                              0&0\cr
                                              0&t_{12}t_{21}
                                             \end{array}\right) \cr
\nonumber
& &-\theta(-k)q^{m+n-1}(1-q^{-2k})r^{(k+1),m,n}\cdot\left(\begin{array}{cc}
                                              t_{12}t_{21}&0\cr
                                              0&0
                                             \end{array}\right)\cr
r^{k,m,n}\cdot\left(\begin{array}{cc}
                                t_{33}&t_{34}\cr
                                t_{43}&t_{44}
                                \end{array}\right)&=&
r^{k,m,n}\cdot\left(\begin{array}{cc}
                                q^{-(m+n)}t_{22}&-q^{-k}t_{21}\cr
                                -q^{-k}t_{12}&q^{m+n}t_{11}
                                \end{array}\right)          \cr
\nonumber
& &-\theta(k)q(1-q^{-2k})r^{(k-1),m,n}\cdot\left(\begin{array}{cc}
                                              t_{12}t_{21}&0\cr
                                              0&0
                                             \end{array}\right) \cr
\nonumber
& &-\theta(-k)q^{m+n-1}(1-q^{-2k})r^{(k+1),m,n}\cdot\left(\begin{array}{cc}
                                              0&0\cr
                                              0&t_{12}t_{21}
                                             \end{array}\right)~~~,
\end{eqnarray}
where $\theta(k)=1$ for $k\geq0$ and $\theta(k)=0$ for $k<0$.
\end{lemma}
\medskip
{\it Proof}. They are obtained by using the following relations
valid for $i<k,j<l$
\begin{equation}
\label{hass}
\begin{array}{rcl}
t_{ij}^n t_{kl}&=&t_{kl}t_{ij}^n -
q^{-1}(1-q^{2n})\,t_{il}t_{kj}t_{ij}^{n-1} \ ,\cr
t_{ij}t_{kl}^m&=&t_{kl}^mt_{ij} + q\,
(1-q^{-2m})\,t_{il}t_{kj}t_{kl}^{m-1}\ . ~~~~~~~~\fidi
\end{array}
\end{equation}

\medskip

By construction
\begin{equation*}
\label{Ccas} \Delta_r := (\id\otimes r)\Delta \ :\
\A(\sphere^7_q)\to\A(\sphere^7_q)\otimes \A(SU(2)_q)
\end{equation*}
defines a right $\A(SU(2)_q)$--coaction on $\A(\sphere^7_q)$. It
satisfies $\Delta_r(uv)=\Delta_r(u)\Delta(v)$, for $u,v\in
\A(\sphere^7_q)$.
The space of functions on the quantum four--sphere $\sph4$ is the
space of coinvariants with respect to this coaction, {\it i.e.}
$\A(\sph4)=\{a\in \A(\sphere^7_q) \, | \, \Delta_r(a) = a\otimes
r(1) \}$. The algebra $\A(\sph4)$ is generated by
$\{a,a^*,b,b^*,R\}$, where
$$a = z_1 z_4^*-z_2 z_3^*,\  b = z_1z_3 + q^{-1}z_2 z_4,
\ R = z_1 z_1^* + z_2 z_2^* \;.$$
They satisfy the following relations
$$
R a = q^{-2} a R \ ,\quad R b = q^2 b R \ ,\quad a b = q^3 b a \ ,\quad
ab^*= q^{-1} b^* a, $$
$$aa^*+q^2bb^* = R(1-q^2R),$$
$$ aa^*= q^2a^*a + (1-q^2)R^2 \ ,\qquad b^*b =
q^4bb^* + (1-q^2)R\;.
$$


\section{The $SU_q(2)$ principal bundle $\sphere^7_q$ over $\sph4$}
In this section we investigate if the structure introduced above
forms a coalgebra-Galois extension, which is essential to define
an (algebraic) quantum principal bundle.

Let $C$ be a coalgebra, $P$ a right $C$-comodule algebra with the
multiplication $m_P: P\ot P\rightarrow P$ and coaction $\dr :
P\rightarrow P\ot C$ .
Let $B\inc P $ be the subalgebra of coinvariants, {\it i.e.}
$B=\{b\in P \,|\, \dr(bp)=b\dr(p), \, \forall\, p\in P\}$.
The canonical left $P$-linear right $C$-colinear map $\chi$
is defined by
\begin{equation*}
\chi:=(m\sb P\ot id)\circ (id\ot\sb B \dr )\, :\; P\ot\sb B P\lra
P\ot C \ , ~p'\ot_B p \mapsto p' \dr p.
\end{equation*}
If $\chi$ is bijective one says that these data form a
coalgebra--Galois extension (see \cite{BH1}).
In this case the translation map is defined as
\begin{equation}
\label{tra} \vartheta : C\rightarrow P\ot_B P \ , \quad
\vartheta (h) = \sum_{[h]}h\o \ot_B h\t := \chi^{-1} (1\ot h)\ .
\end{equation}
If $C$ is a Hopf algebra and $\Delta_R$ is an algebra homomorphism,
the above structure is called a Hopf--Galois extension.
In this case the translation map $\vartheta$ is always
determined by its values on the algebra generators of $C$;
in fact one has that
\begin{equation}\label{t3b}
\vartheta (h  h' ) = \sum_{[h][h']} {h'}\o h\o \ot_B h\t {h'}\t \ .
\end{equation}

We specify now the above framework to our case of
$P=\A(\sphere^7_q)$, $B=\A(\Sigma^4_q)$, $C=\A(SU_q(2))$ and
$\Delta_R=\Delta_r$. Since $\Delta_r$ is not an algebra
homomorphism we are in the more general setting of coalgebra
extensions.
To answer the question of bijectivity of $\chi$ we
cannot simply define the translation map on generators and then
use formula (\ref{t3b}) to extend it to the whole $C$. We shall
instead generalize (\ref{t3b}) by employing  the
$\A(U_q(4))$--module structure defined in (\ref{Ums}). By
considering the coaction $\Delta$ of $\A(U_q(4))$ on
$\A(\sphere^7_q)$, we define the following right
$\A(\sphere^7_q)$-module structure on $\A(\sphere^7_q)\otimes
\A(SU_q(2))$:
$$
(u\otimes x) \triangleleft v = u v_{(0)}\otimes x\cdot v_{(1)}
~~~~u,v\in \A(\sphere^7_q)\;, ~~ x \in \A(SU_q(2))\;,
$$
where $\Delta(v)=\sum_{(v)} v_{(0)}\otimes v_{(1)}$.

By direct computation one can verify that $\chi$ is in fact also right
$\A(\sphere^7_q)$-linear with respect to $\triangleleft$.

The following Proposition is the main result of the paper.

\medskip
\begin{proposition}
\label{s4qq}
Let $C = \A(SU_q(2))\simeq \A(U_q(4))/\R$, $P =  \A(\sphere^7_q)$,
$B = \A(\sph4)$ and $\dr = \Delta_r$ be defined as above. Then the
canonical map $\chi$ is bijective.
\end{proposition}
{\it Proof}. First we prove the surjectivity of $\chi$ by giving a
right inverse, {\it i.e.} a map $\tau:P\otimes C \rightarrow
P\otimes_B P$ such that $\chi\circ\tau=\id$. We will define it on
$1\ot r^{k,m,n}$ and then extend it by left $P$-linearity on $P\ot
C$.

We give an iterative definition. Let $\tau(1\otimes
r(1))=1\otimes_B 1$.
Next assume that for any $k, m, n$ such that $|k|+m+n\leq N $,
$\tau_{k,m,n} = \tau(1\otimes r^{k,m,n})$ is defined so that
$\chi(\tau_{k,m,n})=1\otimes r^{k,m,n}$.
We claim that (in partial analogy to (\ref{t3b})) we can define
$\tau$ on all elements of the basis of degree $N+1$ by the following
formulas :
\begin{equation}\label{t}
\begin{array}{lll}
\tau_{k+1,m,n}  =& q^{2+m+n} z_1^* \tau_{k,m,n} z_1 + q^{2+m+n}
z_2\tau_{k,m,n} z_2^*
& \cr
&+ q^{2} z_3\tau_{k,m,n} z_3^* + z_4^*\tau_{k,m,n} z_4
&\mbox{for\ $k\geq 0$\ ,}\cr
&&\cr
\tau_{k-1,m,n} =& q^4 z_1\tau_{k,m,n} z_1^* + q^2
z_2^*\tau_{k,m,n} z_2 &\cr & + q^{m+n}z_3^*\tau_{k,m,n} z_3 +
q^{m+n}z_4\tau_{k,m,n} z_4^*
&\mbox{for\ $k\leq 0$\ ,}\cr
&&\cr
\tau_{k,m+1,n} =& q^{2+\frac{|k|-k}{2}} z_1^*\tau_{k,m,n} z_2 -
q^{3+\frac{|k|-k}{2}} z_2\tau_{k,m,n} z_1^*
&\cr &
+ q^{1+\frac{|k|+k}{2}} z_3\tau_{k,m,n} z_4^* -
q^{\frac{|k|+k}{2}} z_4^*\tau_{k,m,n} z_3
&\mbox{for\ $k\in\Z$\ ,}\cr
&&\cr
\tau_{k,m,n+1} =& -q^{3+\frac{|k|-k}{2}} z_1\tau_{k,m,n} z_2^* +
q^{2+\frac{|k|-k}{2}} z_2^*\tau_{k,m,n} z_1&\cr &
 - q^{\frac{k+|k|}{2}}z_3^*\tau_{k,m,n} z_4 +
 q^{1+\frac{k+|k|}{2}} z_4\tau_{k,m,n} z_3^*
&\mbox{for\ $k\in\Z$\ .}
\end{array}
\end{equation}
The proof of the good definition of $\tau$ is postponed to Section 4.
Let us apply $\chi$ to r.h.s. of the first equality of (\ref{t}).
The use of the  left and right $P$-linearity of $\chi$ yields
\begin{equation*}
\begin{array}{ll}
&\sum_{j=1}^4 \left(q^{2+m+n} z_1^* z_j \ot r^{k,m,n}\cdot
t_{j1} + q^{2+m+n} z_2 z_j^* \ot r^{k,m,n}\cdot
t_{j2}^*\right. +\cr
&\left.
q^{2} z_3 z_j^* \ot r^{k,m,n}\cdot  t_{j3}^* +
      z_4^* z_j \ot r^{k,m,n}\cdot  t_{j4}\right) =\cr
&
q^{m+n}(q^2 z_1^* z_1 + q^2 z_2 z_2^* + q^2 z_3 z_3^* + z_4^* z_4)
\ot r^{k,m,n}\cdot t_{11} + \cr & (q^{2+m+n}  z_1^*
z_2 - q^{3+m+n} z_2 z_1^* + q^{1-k} z_3 z_4^* - q^{-k} z_4^* z_3)
\ot r^{k,m,n}\cdot  t_{21} =\cr
&q^{m+n}\ot r^{k,m,n}\cdot t_{11} = 1\ot r^{k+1,m,n}\ ,
\end{array}
\end{equation*}
where in the penultimate equality we have used (\ref{s7}).\\
Similarly, $\chi$ applied to r.h.s. of the second equality of
(\ref{t}) yields ($k<0$)
\begin{equation*}
\begin{array}{ll}
&
(q^4 z_1^* z_1 + q^2 z_2^* z_2 + z_3^* z_3 + z_4 z_4^*)
\ot r^{k,m,n}\cdot t_{22} +
\cr &
(-q^{3}  z_1 z_2^* + q^{2} z_2^* z_1 - q^{m+n-k} z_3^* z_4 +
q^{1+m+n-k} z_4 z_3^*)\ot r^{k,m,n}\cdot t_{12} =\cr
&1\ot r^{k,m,n}\cdot t_{22}= 1\ot r^{k-1,m,n} \ .
\end{array}
\end{equation*}
Next, the r.h.s. of the third equality of (\ref{t}), after
application of $\chi$ yields
\begin{equation*}
\begin{array}{ll}
&
q^{\frac{|k|-k}{2}}(q^2 z_1^* z_1 + q^2 z_2 z_2^* + q^2 z_3 z_3^*
+ z_4^* z_4) \ot r^{k,m,n}\cdot t_{12} + \cr &
[q^{2+\frac{|k|-k}{2}} (z_1^* z_2 - q z_2 z_1^*) +
q^{\frac{k+|k|}{2}-m-n} (qz_3 z_4^* -  z_4^* z_3)] \ot
r^{k,m,n}\cdot  t_{22} + \cr &
\theta(k)q^{1+\frac{k+|k|}{2}}(q^{-2k} -1) (q z_3z_4^* -z_4^*z_3)
\ot r^{k-1,m,n}\cdot t_{12}t_{21} =\cr
&1\ot r^{k,m+1,n} \ .
\end{array}
\end{equation*}
Finally, $\chi$ applied to r.h.s. of the fourth equality of
(\ref{t}) yields
\begin{equation*}
\begin{array}{ll}
&
q^{\frac{|k|-k}{2}}(q^4 z_1^* z_1 + q^2 z_2^* z_2 + z_3^* z_3 +
z_4 z_4^*) \ot r^{k,m,n}\cdot t_{21} + \cr &
q^{\frac{|k|-k}{2}}(-q^{3}  z_1 z_2^* + q^{2} z_2^* z_1 -
q^{k+m+n} z_3^* z_4 + q^{1+k+m+n} z_4 z_3^*) \ot r^{k,m,n}\cdot
t_{11} +\cr
&\theta(-k)q^{\frac{|k|+k}{2}}q^{m+n}(1-q^{-2k})(q^{-1}z_3^*z_4-z_4z_3^*) \ot
r^{k+1,m,n}\cdot t_{12}t_{21}= \cr
&1\ot r^{k,m,n+1} \;.
\end{array}
\end{equation*}
The map $\tau$ is then defined by giving its action on the basis
and it satisfies $\chi\circ \tau= {\rm id}$ so that $\chi$ is
surjective.

Injectivity of $\chi$ is a consequence of the $P$--linearity of
the map $\tau$, whose proof is postponed to the next Section. In
fact we have that $\tau(\chi(a\otimes b))=\tau(a\chi(1\otimes
1)\triangleleft b)= a\tau(1\otimes r(1))b=a\otimes b$, so that
$\tau=\chi^{-1}$. \fidi

\bigskip
\bigskip


\section{Proof of good definition and $P$--linearity of $\tau$}

The relations (\ref{t}) provide an iterative definition of the map $\tau$ once
we prove that $\tau_{k,m,n}$ is uniquely defined starting from
$\tau_{k-1,m,n}\;$, $\tau_{k,m-1,n}\;$, $\tau_{k,m,n-1}\;$ .
We will prove the good definition of $\tau$ on the linear basis $r^{k,m,n}$
by induction on $|k|+m+n$. The $P$-linearity will come as an easy consequence
of the induction procedure. Since computations are quite heavy and long we will limit ourselves to
sketch the proof.

We suppose that $\tau(1\otimes r^{k,m,n})=\tau_{k,m,n}$ is well defined for
$|k|+m+n \leq N$ and that the following relations are true for $|k|+m+n\leq
N-1$:
\begin{equation}
\label{facili}\left.
\begin{array}{c}
q^{-(m+n)}z_1\tau_{k+1,m,n}+z_2\tau_{k,m,n+1}=\tau_{k,m,n}z_1\
,\cr -q^{-k}z_3\tau_{k,m,n+1}+z_4\tau_{k+1,m,n}=\tau_{k,m,n}z_4\
,\cr
-qz_1^*\tau_{k,m,n+1}+z_2^*q^{-(m+n)}\tau_{k+1,m,n}=\tau_{k,m,n}z_2^*\
,\cr
z_3^*\tau_{k+1,m,n}+z_4^*q^{-(1+k)}\tau_{k,m,n+1}=\tau_{k,m,n}
z_3^*\ ,
\end{array} \right. k\geq0
\end{equation}

\begin{equation}
\label{facilineg} \left.\begin{array}{c}
q^{-(m+n)}z_3\tau_{k-1,m,n}-z_4\tau_{k,m+1,n}=\tau_{k,m,n} z_3\
,\cr q^{k}z_1\tau_{k,m+1,n}+z_2\tau_{k-1,m,n}=\tau_{k,m,n}z_2\
,\cr
qz_3^*\tau_{k,m+1,n}+z_4^*q^{-(m+n)}\tau_{k-1,m,n}=\tau_{k,m,n}z_4^*\
,\cr
z_1^*\tau_{k-1,m,n}-z_2^*q^{-(1-k)}\tau_{k,m+1,n}=\tau_{k,m,n}
z_1^*\ ,\end{array}\right. k\leq 0
\end{equation}

and the following ones for $|k|+m+n\leq N$:
\begin{equation}
\label{difficili_bis}\left. \begin{array}{c}
q^{m+n}z_2\tau_{k-1,m,n} + q^{m+n+1} \tau_{k-1,m+1,n} z_1 =
\tau_{k,m,n} z_2 \cr z_3 \tau_{k-1,m,n} - q^{-k} \tau_{k-1,m+1,n}
z_4 = \tau_{k,m,n} z_3\cr q^{m+n}z_1^*\tau_{k-1,m,n}-
q^{m+n}\tau_{k-1,m+1,n}z_2^* = \tau_{k,m,n}z_1^*\cr
z_4^*\tau_{k-1,m,n} + q^{1-k}\tau_{k-1,m+1,n} z_3^* =
\tau_{k,m,n}z_4^*\ ,\end{array}\right. k\geq 1
\end{equation}

\begin{equation}
\label{difficilineg_bis}\left. \begin{array}{c}
q^{m+n}z_4\tau_{k+1,m,n} - q^{m+n+1} \tau_{k+1,m,n+1} z_3 =
\tau_{k,m,n} z_4 \cr z_1 \tau_{k+1,m,n} + q^k \tau_{k+1,m,n+1}
z_2 = \tau_{k,m,n} z_1\cr q^{m+n}z_3^*\tau_{k+1,m,n}+
q^{m+n}\tau_{k+1,m,n+1}z_4^* = \tau_{k,m,n}z_3^*\cr
z_2^*\tau_{k+1,m,n} - q^{1+k}\tau_{k+1,m,n+1} z_1^* =
\tau_{k,m,n}z_2^*\ .\end{array}\right. k\leq -1
\end{equation}

It is a straightforward computation to verify the induction hypothesis
for $N=1$ and that (\ref{facili}), (\ref{facilineg}), (\ref{difficili_bis}),
(\ref{difficilineg_bis}) imply the good definition of $\tau_{k,m,n}$ for
$|k|+m+n=N+1$.

In the following subsections we sketch the proof of (\ref{facili}) and
(\ref{difficili_bis}). The set of equations (\ref{facilineg}) and
(\ref{difficilineg_bis}) are proven along the same lines.

The left $P$--linearity of the map $\tau$ is clear by construction.
The right $P$--linearity is a direct consequence of equations
(\ref{facili}), (\ref{facilineg}), (\ref{difficili_bis}),
(\ref{difficilineg_bis}) as can be seen by writing explicitly
$\tau((1\otimes r^{k,m,n})\triangleleft z_\ell )=
\tau(1\otimes r^{k,m,n})z_\ell$ for each $\ell$ and the
corresponding one for $z_\ell^*$.

For each $n\in\Z$ let us define $a_n= z_1z_4^*-q^nz_2z_3^*$ and
$b_n=z_1z_3+q^{n-1}z_2z_4$.

\subsection{Proof of (\ref{facili})}

The following Lemma is preliminary to get the result.

\begin{lemma}
\label{lemma_facili} Relations (\ref{facili}) are true if and only
if the following relations hold for $k\geq0$:
\begin{eqnarray}\label{facili_bis}
(1-q^2 R) \tau_{k,m,n} z_1 = q^{2-m-n} b_{m+n+k} \tau_{k,m,n}z_3^*
+ q^{-m-n} a_{k+m+n} \tau_{k,m,n} z_4  \;, \cr R\tau_{k,m,n}z_4 =
q^{2-k}b_{k+m+n}\tau_{k,m,n}
z_2^*+q^{2+m+n}a^*_{-(m+n+k)}\tau_{k,m,n}z_1\cr
(1-q^2R)\tau_{k,m,n}z_2^* = q^k b^*_{-(m+n+k)}\tau_{k,m,n} z_4 -
q^{3+k} a^*_{-(m+n+k)} \tau_{k,m,n} z_3^*\cr q^4 R
\tau_{k,m,n}z_3^* = q^{2+m+n} b^*_{-(m+n+k)}\tau_{k,m,n}z_1-
q^{3-k} a_{k+m+n}\tau_{k,m,n}z_2^*\;.
\end{eqnarray}
\end{lemma}
{\it Proof.} In order to go from (\ref{facili}) to
(\ref{facili_bis}), substitute in each of (\ref{facili}) the relation
(\ref{t}). All the steps can be retraced back to go in the opposite direction.
\fidi

\medskip

\begin{lemma}
As a consequence of the induction hypothesis, for each $(k,m,n)$
such that $k+m+n=N-1$ we have:
\begin{eqnarray}
\label{facili_tris} (1-q^2R)\tau_{k+1,m,n} = z_4^* \tau_{k,m,n}z_4
+ q^2 z_3\tau_{k,m,n}z_3^*\cr b_{k+m+n+1}\tau_{k+1,m,n} q^{-m-n} =
q z_3 \tau_{k,m,n} z_1 + q^k z_2 \tau_{k,m,n} z_4 \cr a_{k+m+n+1}
\tau_{k+1,m,n} q^{-m-n} = q^{-1}z_4^*\tau_{k,m,n}z_1 - q^{1+k} z_2
\tau_{k,m,n} z_3^*\cr  a^*_{-(m+n+k+1)} \tau_{k+1,m,n}=
q^{-1}z_1^*\tau_{k,m,n}z_4 - q^{-k-1} z_3\tau_{k,m,n} z_2^*\cr
R\tau_{k+1,m,n} = q^{m+n}
(z_1^*\tau_{k,m,n}z_1+z_2\tau_{k,m,n}z_2^*)\cr b^*_{-(m+n+k+1)}
\tau_{k+1,m,n} = q z_1^*\tau_{k,m,n} z_3^* + q^{-k-2}
z_4^*\tau_{k,m,n}z_2^* \;.
\end{eqnarray}
\end{lemma}
{\it Proof.} For instance, in order to get the first one, multiply
on the left by $q^2z_3$ and by $z_4^*$ the fourth and the second
equation of (\ref{facili}), respectively; then add and collect the
terms. All the other relations are obtained in a similar way.
\fidi

Let $(k,m,n)$ be such that $k+m+n=N$ and apply (\ref{facili_tris})
for $(k-1,m,n)$. For example, in order to prove the first of
(\ref{facili_bis}),
\begin{eqnarray*}
(1-q^2R) \tau_{k,m,n} z_1 &=& (z_4^*\tau_{k-1,m,n}z_4 +
q^2z_3\tau_{k-1,m,n}z_3^*)z_1\cr
&=&
(q^{-1}z_4^*\tau_{k-1,m,n}z_1-q^kz_2\tau_{k-1,m,n}z_3^*)z_4+\cr
& & ~~~q^2(q^{k-1}z_2\tau_{k-1,m,n}z_4 + q z_3\tau_{k-1,m,n}z_1)z_3^*\cr
&=& q^{-m-n}(a_{k+m+n}\tau_{k,m,n}z_4 + q^2
b_{k+m+n}\tau_{k,m,n})z_3^*~,
\end{eqnarray*}
where the first and third lines are obtained by using
(\ref{facili_tris}) for $(k-1,m,n)$. All the other equations
(\ref{facili_bis}) are obtained in a similar way. Using Lemma
\ref{lemma_facili} the (\ref{facili}) are then proved.

\subsection{Proof of (\ref{difficili_bis})}

In order to prove them we use the strategy adopted for
(\ref{facili}). By substituting in (\ref{difficili_bis}) the
iterative definition of $\tau$ we get the following Lemma.

\begin{lemma}
Relations (\ref{difficili_bis}) are true if and only if the
following equations hold for $k\geq1$
\begin{eqnarray}
z_2 \tau_{k-1,m,n} (1-q^2R) = - q^{k+2}z_3\tau_{k-1,m,n}
a_{-(k+m+n-1)} + q^{k-1} z_4^*\tau_{k-1,m,n} b_{-(k+m+n-1)}\cr q^2
z_3 \tau_{k-1,m,n}R = q^{m+n} z_1^*
\tau_{k-1,m,n}b_{-(m+n+k-1)}-q^{2-k} z_2 \tau_{k-1,m,n}
a^*_{m+n+k-1}\cr q^{m+n} z_1^* \tau_{k-1,m,n}(1-q^2R) = q^2 z_3
\tau_{k-1,m,n} b^*_{k+m+n-1}+ z_4^*\tau_{k-1,m,n} a^*_{m+n+k-1}\cr
z_4^* \tau_{k-1,m,n} R = q^{2+m+n} z_1^*\tau_{k-1,m,n}
a_{-(m+n+k-1)} + q^{3-k} z_2 \tau_{k-1,m,n} b^*_{m+n+k-1} \;.
\end{eqnarray}\fidi
\end{lemma}

We have the following Lemma:

\begin{lemma}
As a consequence of the induction hypothesis, for each $k+m+n=N$
we have:
\begin{eqnarray}
\label{difficili_tris} \tau_{k,m,n} (1-q^2R) =
z_4^*\tau_{k-1,m,n}z_4 + q^2 z_3
\tau_{k-1,m,n}z_3^*\cr
\tau_{k,m,n}R =
q^{m+n}(z_2\tau_{k-1,m,n}z_2^*+z_1^*\tau_{k-1,m,n} z_1)\cr
\tau_{k,m,n} a_{-(k+m+n)} = q^{-1} z_4^*
\tau_{k-1,m,n}z_1-q^{-k}z_2\tau_{k-1,m,n}z_3^*\cr
\tau_{k,m,n}b_{-(k+m+n)} = q z_3\tau_{k-1,m,n}z_1 +
q^{-k-1}z_2\tau_{k-1,m,n}z_4\cr \tau_{k,m,n}
a^*_{k+m+n}=q^{m+n-1}z_1^*\tau_{k-1,m,n}z_4 -
q^{k+m+n}z_3\tau_{k-1,m,n} z_2^*\cr \tau_{k,m,n}b^*_{k+m+n}=
q^{m+n+1}z_1^*\tau_{k-1,m,n}z_3^*+q^{k+m+n-1}
z_4^*\tau_{k-1,m,n}z_2^*\;.
\end{eqnarray}
\end{lemma}
{\it Proof}. For example, in order to get the first relation,
multiply on the right by $z_2^*$ and by $z_1$ the first and the
third of (\ref{difficili_bis}), respectively, and then add. The
other relations are shown in a similar way. \fidi

Let $(k,m,n)$ be such that $k+m+n=N+1$ and let us apply (\ref{difficili_tris})
for $(k-1,m,n)$.
For example, in order to show the first of (\ref{difficili_tris})
\begin{eqnarray*}
z_2 \tau_{k-1,m,n}(1-q^2R) &=&z_2z_4^*\tau_{k-2,m,n}z_4 + q^2z_2
z_3\tau_{k-2,m,n} z_3^*\cr
&=& q^{k-1} z_4^* (q^{-k}z_2\tau_{k-2,m,n}z_4 + q z_3
\tau_{k-2,m,n}z_1) -\cr
& & ~~ q^{k+2} z_3 (-q^{1-k} z_2\tau_{k-2,m,n} z_3^* + q^{-1}z_4^*
\tau_{k-2,m,n}z_1 ) \cr
&=& q^{k-1} z_4^* \tau_{k-1,m,n} b_{-(k+m+n-1)} - q^{k+2} z_3
\tau_{k-1,m,n} a_{-(k-1+m+n)} \;,
\end{eqnarray*}
where we used the first, the third and the fourth of the
(\ref{difficili_tris}) for $(k-1,m,n)$. All the others relations
are obtained in a similar way. This concludes the proof of
(\ref{difficili_bis}).


\section{Conclusions}
The bijectivity of the canonical map is a key result to study the noncommutative
geometry of this fibration.

To this purpose it is relevant the notion of \emph{principal Galois extension}
introduced in \cite{BH2,HMS} which consists of a coalgebra Galois
extension with the following additional requirements (the notations are the
same as in the beginning of Section 3): {\it i\,}) the entwining map $\psi :
C\ot P \rightarrow P \ot C$, $c \ot p\rightarrow \chi(\chi^{-1}(1\ot c)\, p)$
is bijective; {\it ii\,}) there exists a group like element $e \in C$ such that
$\Delta_R (p)=\psi(e \ot p)$; {\it iii\,}) there exists a strong connection.

In our case the points {\it i}) and {\it ii}) are easily satisfied. In fact
$\psi(c \ot p) =(1 \ot c)\triangleleft p$ can be inverted by $\psi^{-1}(p \ot
c)=c\cdot S^{-1}(p_{(1)})\ot p_{(0)}$ and $e= r(1)$. It is crucial
that $P=\A(\sphere^7_q)$ is a $\A(U_q(4))$--comodule  and $C=\A(SU_q(2))$
a $\A(U_q(4))$--module.

A result (Theorem 2.44) contained in the paper in preparation \cite{schsch}
solves the point {\it iii\,}): if $C$ is a cosemisimple coalgebra, right
$H$--module quotient of an Hopf algebra $H$ with bijective antipode then the
surjectivity of the canonical map implies its injectivity and the existence
of a strong connection (or equivariant projectivity in their terminology).

The translation map $\vartheta$ can be lifted to $\ell: \A(SU_q(2)) \rightarrow
\A(\sphere^7_q) \ot \A(\sphere^7_q)$ by making use of the relations (\ref{t}).
We conjecture that the lifted map $\ell$ is bicolinear with respect
to $(1 \ot \Delta_r)$ and $(\bar{\Delta}_r \ot 1)$ where $\bar{\Delta}_r:p
\rightarrow r(S^{-1}(p_{(1)}))\ot p_{(0)}$ \cite{HMS}, so that it gives an
explicit expression to the strong connection.

A relevant consequence of principality is that it makes possible, in the
quantum setting, the construction of associated vector bundles.
Given, in fact, any finite dimensional left corepresentation
of $C$, ($\rho:V\to C \ot V$), the cotensor product $P\Box_\rho V$ turns out
to be a finitely generated and projective  $B$--module.

In our case for any corepresentation of spin $j$ of $\A(SU_q(2))$ one has a
\emph{quantum vector bundle} and the corresponding class in the positive cone
of $K_0(\A(\sph4))$.

In \cite{BCT1} it was already computed the vector bundle associated to the
spin $1/2$ corepresentation and its pairing with the $K$--homology
generators. The non triviality of this paring allows to conclude that the
coalgebra Galois extension is non--cleft. Using the formula for the
Chern--Connes character given in \cite{BH2} we have all we need to compute
the Chern numbers of the generic associated bundles. This will be the content
of a future work.

\bigskip
\bigskip

{\bf Acknowledgements.} We want to thank P.M. Hajac for enlightening
discussions, for explaining us his results in \cite{BH2} and for providing us
useful bibliographical informations, and P. Schauenburg and H.--J. Schneider
for anticipating us the results contained in \cite{schsch}.

\end{document}